\documentclass[12pt]{ert-l}
\usepackage{amssymb,amscd}
\usepackage[margin=2.8cm]{geometry}
\newtheorem{thm}[subsection]{Theorem}
\newtheorem{lem}[subsection]{Lemma}
\newtheorem{prop}[subsection]{Proposition}
\newtheorem{cor}[subsection]{Corollary}
\theoremstyle{definition}

\theoremstyle{remark}

\numberwithin{equation}{subsection}

\newcommand{\N}{{\mathbb N}}
\newcommand{\Z}{{\mathbb Z}}
\newcommand{\Q}{{\mathbb Q}}

\newcommand{\ad}{\operatorname{ad}}
\newcommand{\End}{\operatorname{End}}

\newcommand{\GL}{{\sf GL}}

\newcommand{\gl}{\mathfrak{gl}}

\newcommand{\g}{\mathfrak{g}}
\newcommand{\h}{\mathfrak{h}}
\newcommand{\divided}[2]{#1^{(#2)}}
\newcommand{\sqbinom}[2]{\begin{bmatrix}#1\\#2\end{bmatrix}}
\newcommand{\U}{\mathbf{U}}
\let\sect=\S
\renewcommand{\S}{\mathbf{S}}
\newcommand{\T}{\mathbf{T}}
\newcommand{\B}{\mathbf{B}}
\newcommand{\A}{\mathcal{A}}
\renewcommand{\i}{{1}}
\renewcommand{\b}{{\sf b}}

\newcommand{\bil}[2]{\langle #1, #2 \rangle}
\newcommand{\opp}{{\operatorname{opp}}}
\newcommand{\lsub}[2]{{}_{#2}{#1}}

\begin{document}
\title{Presenting generalized $q$-Schur algebras}
\author{Stephen Doty}
\address{Loyola University Chicago, Chicago, Illinois 60626 U.S.A.}
\email{doty@math.luc.edu}
\date{August 31, 2002}
\subjclass{17B37, 16W35, 81R50} 
\keywords{Schur algebras, q-Schur algebras, generalized Schur algebras, 
quantized enveloping algebras}

\begin{abstract} 
We obtain a presentation by generators and relations for generalized
Schur algebras and their quantizations. This extends earlier results
obtained in the type $A$ case.  The presentation is compatible with
Lusztig's modified form ${\bf\dot{\U}}$ of a quantized enveloping
algebra. We show that generalized Schur algebras inherit a canonical
basis from ${\bf\dot{\U}}$, that this gives them a cellular structure, and
thus they are quasihereditary over a field. 
\end{abstract}
\maketitle

\parskip=2pt
\allowdisplaybreaks

\section*{Introduction}

In \cite{Donkin:SA1} Donkin defined the notion of a generalized Schur
algebra for an algebraic group, depending on the group and a finite
saturated subset $\pi$ of dominant weights. He also showed how to
construct the generalized Schur algebra from the enveloping algebra of
the complex Lie algebra of the same type as the given algebraic group,
by an appropriate modification of the construction of Chevalley groups. 

The purpose of this paper is to give a presentation by generators and
relations for generalized Schur algebras and their quantizations.  The
presentation has the same form as \cite[Theorems 1.4, 2.4]{DG:PSA},
for Schur algebras and $q$-Schur algebras in type $A$.

We approach this problem the other way around. First we define (in
\sect\ref{Spi}) an algebra $\S(\pi)$ (over $\Q(v)$, $v$ an
indeterminate) by generators and relations.  It depends only on a
Cartan matrix (of finite type) and a given saturated set $\pi$ of
dominant weights. We prove that this algebra is a finite-dimensional
semisimple quotient of the quantized enveloping algebra $\U$
determined by the Cartan matrix, and that it is a $q$-analogue of a
generalized Schur algebra in Donkin's sense.  We show that $\S(\pi)$
is isomorphic with the algebra ${\bf\dot{\U}}/ {\bf\dot{\U}}[P]$ constructed by
Lusztig in \cite{Lusztig:book}, where $P$ is the complement of $\pi$
in the set of dominant weights and where ${\bf\dot{\U}}$ is his ``modified
form'' of $\U$, using his ``refined Peter-Weyl theorem,'' which for
convenience we recall in \sect\ref{recollect}.  It follows that
$\S(\pi)$ inherits a canonical basis and a cell datum (in the sense of
Graham and Lehrer \cite{GL}) from ${\bf\dot{\U}}$. This provides $\S(\pi)$
with an extremely rigid structure, which essentially determines its
representation theory in all possible specializations. In particular,
$\lsub{\S}{R}(\pi)$ is quasihereditary over any field $R$.  These
results are contained within \sect\sect\ref{GqSA}--\ref{basechange}.

In \sect\ref{classical} we consider the classical case. First we
define an algebra $S(\pi)$ (over $\Q$) by generators and
relations. The defining presentation of $S(\pi)$ is obtained from the
defining presentation of $\S(\pi)$ by simply setting $v=1$.  Then we
show that $S(\pi)$ is a generalized Schur algebra.  It has essentially
the same cell datum as does $\S(\pi)$, and it is quasihereditary when
specialized to any field. (This also follows from results in
\cite{Donkin:SA1}.)  The arguments are quite similar to those used in
the quantum case, so many of them are omitted or sketched.

In \sect\ref{examples} we consider some natural examples, and in 
\sect\ref{integral} we obtain generators and relations for integral
forms of generalized Schur algebras.

We mostly follow notational conventions of \cite{Lusztig:book}.

\section{The algebra $\S(\pi)$}\label{Spi}

\subsection{} \label{Cartan:matrix}
Let $(a_{ij})_{1\le i,j \le n}$ be a Cartan matrix of finite type. We
do not assume it is indecomposable.  The entries $a_{ij}$ are integers
and there is a vector $(d_1, \dots, d_n)$ with entries $d_i \in
\{1,2,3\}$ such that the matrix $(d_ia_{ij})$ is symmetric and
positive definite, $a_{ii}=2$ for all $i$, and $a_{ij}\le 0$ for all
$i \ne j$.

\subsection{} \label{g:relations}
Let $\g$ be the corresponding finite-dimensional semisimple Lie algebra
over the rational field $\Q$. As we know, $\g$ is the Lie algebra
given by the generators $e_i, f_i, h_i$ ($1\le i \le n$) and relations
\begin{gather}
{}[h_i,h_j]=0, \quad [e_i,f_j]=\delta_{ij}h_i, 
\tag{a} \\ 
[h_i, e_j] = a_{ij}e_j, \quad [h_i, f_j] = -a_{ij}f_j, 
\tag{b} \\ 
(\ad e_i)^{1-a_{ij}} e_j = 0 = (\ad f_i)^{1-a_{ij}} f_j \quad (i \ne j).  
\tag{c} 
\end{gather} 
We fix the Cartan subalgebra $\h$ spanned by the $h_i$. Let $\Phi$ be
the root system of $\g$ with respect to $\h$, and let $W$ be the Weyl
group of $\Phi$. We fix ordered bases $\Pi = \{\alpha_1, \dots, \alpha_n\}
\subset \h^*$, $\Pi^\vee  = \{\alpha^\vee_1, \dots, \alpha^\vee_n\}
\subset \h$ such that $\bil{\alpha_i^\vee}{\alpha_j} = a_{ij}$ 
for all $i,j$.

\subsection{} \label{qfactorial}
Let $v$ be an indeterminate, and set $\A = \Z[v,v^{-1}]$, with
quotient field $\Q(v)$. Set $v_i = v^{d_i}$.  Given $a\in \Z$, $t\in
\N$ set $[a]_i = (v_i^a - v_i^{-a})/(v_i - v_i^{-1})$,
\begin{equation} 
[t]_i^! = \prod_{s=1}^t  \frac{v_i^s - v_i^{-s}}{v_i - v_i^{-1}}, \quad
\sqbinom{a}{t}_i = \frac{\prod_{s=0}^{t-1} (v_i^{a-s} - v_i^{-a+s})}
{\prod_{s=1}^t (v_i^s - v_i^{-s})}.
\end{equation}
The subscript $i$ can be omitted when $d_i=1$. 

\subsection{} \label{weights}
Let $X = \{ \lambda \in \h^* \mid \bil{\alpha_i^\vee}{\lambda} \in \Z,
\text{ all $i = 1, \dots, n$} \}$ (the weight lattice) and $X^+ = \{
\lambda \in X \mid \bil{\alpha_i^\vee}{\lambda} \ge 0 \text{ all $i =
1, \dots, n$} \}$ (the set of dominant weights).  For $\lambda,\mu \in
X$ write $\lambda \le \mu$ if $\mu-\lambda \in \N\alpha_1 + \cdots +
\N\alpha_n$.  This defines a partial order on $X$.

\subsection{} \label{S:relations}
Let $\pi$ be a finite ideal in the poset $X^+$ ($\lambda \le \mu$ for
$\mu \in \pi$, $\lambda\in X^+$ implies $\lambda \in \pi$). Such sets
of weights are sometimes called {\em saturated}.  We define an algebra
(associative with $1$) over $\Q(v)$ given by generators $E_i, F_i$ ($1
\le i \le n$), $\i_\lambda$ ($\lambda \in W\pi$) and relations
\begin{gather}
\i_\lambda \i_\mu = \delta_{\lambda\mu} \i_\lambda, \quad
\sum_{\lambda\in W\pi} \i_\lambda = 1 
\tag{a} \\
E_i F_j - F_j E_i = \delta_{ij} \sum_{\lambda\in W\pi}
[\bil{\alpha_i^\vee}{\lambda}]_i \i_\lambda 
\tag{b} \\
E_i \i_\lambda =
\begin{cases}
\i_{\lambda+\alpha_i} E_i &
   \text{if $\lambda+\alpha_i \in W\pi$}\\
0 & \text{otherwise}
\end{cases} \tag{c} \\
F_i \i_\lambda =
\begin{cases}
\i_{\lambda-\alpha_i} F_i &
   \text{if $\lambda-\alpha_i \in W\pi$}\\
0 & \text{otherwise}
\end{cases} \tag{d} \\
\i_\lambda E_i =
\begin{cases}
E_i \i_{\lambda-\alpha_i} &
   \text{if $\lambda-\alpha_i \in W\pi$}\\
0 & \text{otherwise}
\end{cases} \tag{e} \\
\i_\lambda F_i =
\begin{cases}
F_i \i_{\lambda+\alpha_i} &
   \text{if $\lambda+\alpha_i \in W\pi$}\\
0 & \text{otherwise}
\end{cases} \tag{f} \\
\sum_{s=0}^{1-a_{ij}} (-1)^s \sqbinom{1-a_{ij}}{s}_i 
E_i^{1-a_{ij}-s}E_j E_i^s = 0 \quad (i\ne j)
\tag{g} \\
\sum_{s=0}^{1-a_{ij}} (-1)^s \sqbinom{1-a_{ij}}{s}_i
F_i^{1-a_{ij}-s}F_j F_i^s = 0 \quad (i\ne j).
\tag{h} 
\end{gather}
Denote this algebra by $\S=\S(\pi)$. It depends only on the Cartan
matrix $(a_{ij})$ and the saturated set $\pi$.

\section{Recollections}\label{recollect}

\subsection{} \label{U:relations}
Let $\U$ be the algebra (associative with $1$) given by generators
$E_i$, $F_i$, $K_i$, $K_i^{-1}$ ($1\le i \le n$) and relations 
\begin{gather}
K_i K_j =K_j K_i, \quad K_i K_i^{-1} = 1 = K_i^{-1}K_i \tag{a}\\
E_iF_j - F_jE_i = \delta_{ij}\frac{K_i-K_i^{-1}}{v_i-v_i^{-1}} \tag{b}\\
K_i E_j = v_i^{a_{ij}}E_jK_i, \quad K_i F_j = v_i^{-a_{ij}}F_jK_i \tag{c}\\
\sum_{s=0}^{1-a_{ij}} (-1)^s \sqbinom{1-a_{ij}}{s}_i 
E_i^{1-a_{ij}-s}E_j E_i^s = 0 \quad (i\ne j)\tag{d}\\
\sum_{s=0}^{1-a_{ij}} (-1)^s \sqbinom{1-a_{ij}}{s} _i
F_i^{1-a_{ij}-s}F_j F_i^s = 0 \quad (i\ne j).\tag{e}
\end{gather}
Then $\U$ is the Drinfeld-Jimbo quantized enveloping algebra
corresponding to the Cartan matrix $(a_{ij})$. 

\subsection{} \label{modules}
Let $M$ be a finite-dimensional $\U$-module. Then 
$M = \oplus_{\lambda,\sigma} M^\lambda_\sigma$ for $\lambda\in X$,
$\sigma: \{1,\dots,n\} \to \{1,-1\}$. Here 
$$
M^\lambda_\sigma = \{ m\in M \mid K_i m =
\sigma(i)v_i^{\bil{\alpha_i}{\lambda}} m , i = 1, \dots, n \}
$$
is the $(\lambda,\sigma)$ weight space. We write $M_\sigma =
\oplus_{\lambda\in X} M^\lambda_\sigma$; then $M = \oplus M_\sigma$
where the sum is over all maps $\sigma: \{1,\dots,n\} \to \{1,-1\}$.
One says that $M$ has type $\sigma$ if $M=M_\sigma$. In case
$\sigma(i)=1$ for all $i$ then one says that $M$ has type $\bf1$.

It is known that every finite-dimensional $\U$-module is completely
reducible. Moreover, for each $\sigma$ as above there is an
equivalence of categories between finite-dimensional $\U$-modules of
type $\bf1$ and those of type $\sigma$. Thus one usually confines one's
attention to the type $\bf1$ modules.

Every type $\bf1$ simple $\U$-module is a highest weight module of
highest weight $\lambda$, for some $\lambda \in X^+$.  We denote this
module by $\Lambda_\lambda$.  For $\lambda \ne \lambda'$,
$\Lambda_\lambda$ is not isomorphic to $\Lambda_{\lambda'}$.

\subsection{} \label{Udot}
In Lusztig's book \cite[Part IV]{Lusztig:book}, a ``modified form''
${\bf\dot{\U}}$ of $\U$ was introduced, and it was shown how to extend the
canonical basis from the plus part of $\U$ to a canonical basis
${\bf\dot{\B}} = \sqcup_{\lambda \in X^+} {\bf\dot{\B}}[\lambda]$ on all of
${\bf\dot{\U}}$. (See \cite[29.1.1]{Lusztig:book} for the definition of
${\bf\dot{\B}}[\lambda]$.)  The $\U$-modules of type $\bf1$ having weight
space decompositions can be regarded naturally as modules for
${\bf\dot{\U}}$; on the other hand, more exotic $\U$-modules, such as those
without weight decomposition, cannot be regarded as modules for
${\bf\dot{\U}}$.  The algebra ${\bf\dot{\U}}$ in general does not have a unit
element, but it does have a family $(1_\lambda)_{\lambda\in X}$ of
orthogonal idempotents such that ${\bf\dot{\U}} = \oplus_{\lambda,\lambda'
\in X} 1_{\lambda'} {\bf\dot{\U}} 1_\lambda$.  In a sense, the family
$(1_\lambda)$ serves as a replacement for the identity.  We shall need
the following result from \cite[29.3.3]{Lusztig:book}, which Lusztig
calls the ``refined Peter-Weyl theorem.''

\begin{thm} \label{Peter-Weyl}
Given $\lambda \in X^+$, define ${\bf\dot{\U}}[\ge \lambda]$ (resp,
${\bf\dot{\U}}[>\lambda]$ as the set of all $u\in {\bf\dot{\U}}$ with the
property: if $\lambda' \in X^+$ and $u$ acts on $\Lambda_{\lambda'}$
by a non-zero linear map, then $\lambda' \ge \lambda$ (resp.,
$\lambda' > \lambda$).  Then:

(i) ${\bf\dot{\U}}[\ge \lambda]$ and ${\bf\dot{\U}}[>\lambda]$ are two-sided
ideals of ${\bf\dot{\U}}$, which are generated as vector spaces by their
intersections with ${\bf\dot{\B}}$. The quotient algebra ${\bf\dot{\U}}[\ge
\lambda]/ {\bf\dot{\U}}[>\lambda]$ is isomorphic (via the action of
${\bf\dot{\U}}$ on $\Lambda_\lambda$) to the algebra $\End(\Lambda_\lambda)$. 
Let $p:{\bf\dot{\U}}[\ge \lambda] \to {\bf\dot{\U}}[\ge \lambda] / {\bf\dot{\U}}[>\lambda]$
be the natural projection.

(ii) There is a unique direct sum decomposition of ${\bf\dot{\U}}[\ge
\lambda]/ {\bf\dot{\U}}[>\lambda]$ into a direct sum of simple left
${\bf\dot{\U}}$-modules such that each summand is generated by its
intersection with the basis $p({\bf\dot{\B}}[\lambda])$ of ${\bf\dot{\U}}[\ge
\lambda]/ {\bf\dot{\U}}[>\lambda]$.

(iii) There is a unique direct sum decomposition of ${\bf\dot{\U}}[\ge
\lambda]/ {\bf\dot{\U}}[>\lambda]$ into a direct sum of simple right
${\bf\dot{\U}}$-modules such that each summand is generated by its
intersection with the basis $p({\bf\dot{\B}}[\lambda])$ of ${\bf\dot{\U}}[\ge
\lambda]/ {\bf\dot{\U}}[>\lambda]$.

(iv) Any summand in the decomposition (ii) and any summand in the
decomposition (iii) have an intersection equal to a line consisting of
all multiples of some element in the basis $p({\bf\dot{\B}}[\lambda])$.
This gives a map from the set of pairs consisting of a summand in the
decomposition (ii) and one in the decomposition (iii), to the set
$p({\bf\dot{\B}}[\lambda])$. This map is a bijection.
\end{thm}

\subsection{} \label{cells}
For each $\lambda \in X^+$ let $M(\lambda)$ be any finite set in
bijective correspondence with a basis of weight vectors for
$\Lambda_\lambda$. We could simply take $M(\lambda)$ to be such a basis,
for instance. Then the minimal left ideals in the decomposition (ii)
above are indexed by elements of $M(\lambda)$; the same is true of the
minimal right ideals in the decomposition (iii).  According to part (iv)
of the above theorem, there is a unique element $\b^\lambda_{S,T}$ of
${\bf\dot{\B}}[\lambda]$ corresponding to the intersection of the minimal
left ideal $L_T$ indexed by $T$ in the decomposition (ii) with the
minimal right ideal $R_S$ indexed by $S$ in the decomposition
(iii). Then in terms of this notation we have
\begin{equation}
{\bf\dot{\B}}[\lambda] = \{ \b^\lambda_{S,T} \mid S,T \in M(\lambda \};
\end{equation}
moreover, $L_T$ and $R_S$ are spanned, respectively, by the sets
\begin{equation}
\{\b^\lambda_{S,T} \mid  S \in M(\lambda)\}; \quad
\{\b^\lambda_{S,T} \mid  T \in M(\lambda)\}.
\end{equation} 
These sets (each of cardinality $\dim \Lambda_\lambda$) are the left and
right cells of \cite[29.4.1]{Lusztig:book}. The sets
${\bf\dot{\B}}[\lambda]$ are the two-sided cells.

\subsection{} \label{involution}
We note that by \cite[23.1.6]{Lusztig:book} there are maps $\sigma,
\omega: {\bf\dot{\U}} \to {\bf\dot{\U}}$ (induced from correponding maps $\U \to
\U$) such that $\iota = \sigma\omega = \omega\sigma$ is a
$\Q(v)$-linear anti-involution on ${\bf\dot{\U}}$.  From the proof of
\cite[29.3.3(e)]{Lusztig:book} we see that $\iota(\b^\lambda_{S,T}) =
\b^\lambda_{T,S}$ for all $S,T \in M(\lambda)$; in particular, $\iota$
induces an involution on the vector space ${\bf\dot{\U}}[\ge \lambda]/
{\bf\dot{\U}}[>\lambda]$ which interchanges the minimal left and right
ideals in the decompositions (ii), (iii) in the preceding theorem.

\subsection{} \label{cellular}
For $\lambda\in X^+$ and a fixed $T \in M(\lambda)$ the elements
$p(\b^\lambda_{S,T})$ for $S\in M(\lambda)$ span the ${\bf\dot{\U}}$-module
$L_T$; thus for any $u \in {\bf\dot{\U}}$ we have
\begin{equation}
u\, p(\b^\lambda_{S,T}) = \sum_{S' \in M(\lambda)} 
r_u(S',S)\, p(\b^\lambda_{S',T})
\end{equation}
where $r_u(S',S) \in \Q(v)$ does not depend on $T$. It follows that
in ${\bf\dot{\U}}$ we have
\begin{equation}\label{rdef}
u\, \b^\lambda_{S,T} = \sum_{S' \in M(\lambda)} r_u(S',S)\, \b^\lambda_{S',T} 
\pmod{{\bf\dot{\U}}[>\lambda]}.
\end{equation} 
It follows that the basis ${\bf\dot{\B}} = \{ \b^\lambda_{S,T} \}$ is a
cellular basis for the algebra ${\bf\dot{\U}}$; more precisely, the datum
$(X^+, M, \b, \iota)$ is a cell datum for ${\bf\dot{\U}}$, in the sense of
Graham and Lehrer \cite{GL} (with respect to the opposite partial
order on $X^+$). We note that this cell datum is of profinite type, in
the sense of \cite[2.1.2]{RGreen}; hence the algebra ${\bf\dot{\U}}$ can be
completed to a procellular algebra $\widehat\U$.  The algebra
$\widehat{\U}$ was studied in \cite{BLM} in the type $A$ case.

\subsection{} \label{Aform}
Let $\lsub{{\bf\dot{\U}}}{\A}$ be the set of all finite $\A$-linear
combinations of elements of the canonical basis ${\bf\dot{\B}}$. This is an
$\A$-subalgebra of ${\bf\dot{\U}}$; it is generated by all
$\divided{E_i}{m} 1_\lambda$, $\divided{F_i}{m} 1_\lambda$, for
various $\lambda\in X$, $1\le i \le n$, $m \ge 0$, where
\begin{equation}
\divided{E_i}{m} = E_i^m / ([m]_i^!), 
\divided{F_i}{m} = F_i^m / ([m]_i^!).
\end{equation}
Lusztig \cite[25.4]{Lusztig:book} has shown that, with respect the
canonical basis ${\bf\dot{\B}}$, the structure constants of ${\bf\dot{\U}}$
lie in $\A$. It follows that, for $u \in \lsub{{\bf\dot{\U}}}{\A}$, 
\begin{equation}\label{structure}
\text{the elements $r_u(S',S)$ in \ref{cellular}\eqref{rdef} lie in $\A$}
\end{equation}
for all $S,S' \in M(\lambda)$. Hence, $\lsub{{\bf\dot{\U}}}{\A}$ a cellular
algebra with the same cell datum $(X^+, M, \b, \iota)$ considered above.

\section{Generalized $q$-Schur algebras}\label{GqSA}

We show that the algebra $\S=\S(\pi)$ is (a quantization of) a
generalized Schur algebra in the sense of Donkin \cite{Donkin:SA1}.

\subsection{}\label{K_i}
For convenience, set $\lambda_i = \bil{\alpha_i^\vee}{\lambda}$ for
all $1\le i\le n$, $\lambda\in X$. Note that $\lambda \in X$ is
uniquely determined by its vector $(\lambda_1, \dots, \lambda_n)$ of
values on the coroots $\alpha_i^\vee$. In $\S=\S(\pi)$ we define
elements
\begin{equation}
K_i = \sum_{\lambda\in W\pi} v_i^{\lambda_i}
\i_\lambda;\quad  K_i^{-1} = \sum_{\lambda\in W\pi}
v_i^{-\lambda_i} \i_\lambda
\end{equation}
for any $1 \le i \le n$. Then one verifies immediately from
\ref{S:relations}(a) and the definition that $K_i$, $K_i^{-1}$ satisfy
relation \ref{U:relations}(a); in particular the $K_i$ all commute.

\begin{lem}\label{Kgen}
Let $\S^0$ be the subalgebra of $\S$ generated by $K_1, \dots,
K_n$.

(i) The idempotents $\i_\lambda$ ($\lambda\in  W\pi$) lie
within $\S^0$.

(ii) The $K_i^{-1}$ lie in $\S^0$.

(iii) The $\i_\lambda$ ($\lambda\in W\pi$) form a basis of $\S^0$. 
\end{lem}

\begin{proof}\newcommand{\lp}{\lambda^\prime}
We set $\Gamma(i,\lambda) = \{ \mu\in W\pi \mid \mu_i = \lambda_i
\}$, and $J_i^\lambda = \prod_\mu (K_i - v_i^{\mu_i})$ where the 
product is taken over all $\mu \in W\pi - \Gamma(i,\lambda)$.
We have equalities
\begin{equation}
\begin{aligned} 
J_i^\lambda &= \prod_\mu \left( \sum_{\lp \in W\pi} v_i^{\lp_i}\i_{\lp} -
v_i^{\mu_i} \sum_{\lp \in W\pi} \i_{\lp} \right) \\
&=  \prod_\mu \left( \sum_{\lp \in W\pi} (v_i^{\lp_i} -
v_i^{\mu_i}) \i_{\lp} \right) \\
&= \sum_{\lp \in W\pi} \prod_\mu  (v_i^{\lp_i} -
v_i^{\mu_i}) \i_{\lp}
\end{aligned}
\end{equation}
where all products are taken over $\mu \in W\pi -
\Gamma(i,\lambda)$ and where we have used the idempotent orthogonality
relations \ref{S:relations}(a) to interchange the sum and product.
Noting that the product in the sum on the last line above vanishes for
any $\lp \in W\pi - \Gamma(i,\lambda)$, we obtain the expression
\begin{equation}
J_i^\lambda = \sum_{\lp \in \Gamma(i,\lambda)} \prod_\mu  
(v_i^{\lp_i} - v_i^{\mu_i}) \i_{\lp} 
\end{equation}
where the product in this sum is a nonzero {\em constant}, since
$\lp_i = \lambda_i$ for all $\lp \in \Gamma(i,\lambda)$.  This proves
that $J_i^\lambda$ is (up to a nonzero scalar) the sum of all
idempotents $\i_{\lp}$ for which $\lp_i = \lambda_i$. This property
holds for all $i$. Thus it follows that the product $J_1^\lambda
\cdots J_n^\lambda$ is, up to a nonzero scalar multiple, equal 
to $\i_\lambda$, since $\i_\lambda$ is the unique idempotent appearing
in each of the sums in the product. 

By definition $J_i^\lambda$ belongs to the subalgebra of $\S(\pi)$
generated by $K_i$, so the result of the previous paragraph shows that
$\i_\lambda$ (for any $\lambda \in W\pi$) lies within the subalgebra
of $\S(\pi)$ generated by all $K_1, \dots, K_n$. This proves part (i).

Part (ii) follows from part (i) and the definition of $K_i^{-1}$.

By definition of the $K_i$ we see that the subalgebra of $\S$
generated by the $\i_\lambda$ contains the $K_i$. By part (i) this
subalgebra equals $\S^0$. Part (iii) now follows from the fact that
the $\i_\lambda$ form a family of orthogonal idempotents.
\end{proof}

\subsection{} \label{UtoS}
The algebra $\S(\pi)$ is generated by all $E_i, F_i, K_i, K_i^{-1}$
($1\le i \le n$). In fact, by Lemma \ref{Kgen}(ii), it is generated by
the $E_i, F_i, K_i$.  We have already observed that the generators
satisfy relation \ref{U:relations}(a).  From \ref{S:relations}(b) and
the definitions we obtain equalities
\begin{equation}
\begin{aligned}
E_i F_j - F_j E_i &= \delta_{ij} \sum_{\lambda \in W\pi}
[\bil{\alpha_i^\vee}{\lambda}]_i \i_\lambda \\
&= \delta_{ij} \sum_{\lambda \in W\pi} 
\frac{v_i^{\bil{\alpha_i^\vee}{\lambda}} - v_i^{-\bil{\alpha_i^\vee}{\lambda}}}
{v_i - v_i^{-1}} \,\i_\lambda \\
&= \delta_{ij} \frac{K_i - K_i^{-1}}{v_i - v_i^{-1}},
\end{aligned}
\end{equation}
which shows that the elements $K_i, K_i^{-1}, E_i, F_i$ must satisfy
\ref{U:relations}(b).

We set $\i_\mu = 0$ for any $\mu \in X-W\pi$. This gives a meaning to
the symbol $\i_\lambda$ for all $\lambda \in X$. With this convention
we can write the relations \ref{S:relations}(c)--(f) more compactly,
in the form:
\begin{equation}
E_i \i_\lambda = \i_{\lambda+\alpha_i} E_i, \quad 
F_i \i_\lambda = \i_{\lambda-\alpha_i} F_i
\end{equation}
for all $i$ and all $\lambda\in X$. Then we have
\begin{equation}
\begin{aligned}
K_i E_j &= \sum_{\lambda\in X} v_i^{\bil{\alpha_i^\vee}{\lambda}}
\i_\lambda E_j \\ 
&= \sum_{\lambda\in X}
v_i^{\bil{\alpha_i^\vee}{\lambda}} E_j \i_{\lambda-\alpha_j}\\
&= \sum_{\lambda\in X}
v_i^{\bil{\alpha_i^\vee}{\lambda+\alpha_j}} E_j \i_{\lambda}\\
&= v_i^{\bil{\alpha_i^\vee}{\alpha_j}} E_j \sum_{\lambda\in X}
v_i^{\bil{\alpha_i^\vee}{\lambda}} \i_{\lambda}\\
&= v_i^{a_{ij}} E_j K_i
\end{aligned} 
\end{equation}
which proves the first part of relation \ref{U:relations}(c). The
other part of \ref{U:relations}(c) is proved by the analogous
calculation. So the generators $K_i, K_i^{-1}, E_i, F_i$ satisfy
\ref{U:relations}(c).  They also satisfy relations
\ref{U:relations}(d), (e) since those relations are identical with
\ref{S:relations}(g), (h).  We have proved the following result.

\begin{prop}\label{image}
The algebra $\S(\pi)$ is a homomorphic image of $\U$, via the
homomorphism sending the elements $E_i, F_i, K_i, K_i^{-1}$ of $\U$ to
the corresponding elements of $\S(\pi)$ denoted by the same symbols.
\end{prop}

\subsection{}
Set $p_i(X) = \prod_{\mu\in W\pi} (X - v_i^{\mu_i}) \in \Q(v)[X]$,
where $X$ is a formal indeterminate.

\begin{lem} \label{finiteness}
In the algebra $\S(\pi)$ we have:
\par(i) The $K_i$ satisfy the polynomial identity $p_i(K_i) = 0$
for $i = 1, \dots, n$. 
\par(ii) The $E_i$, $F_i$ are nilpotent.
\end{lem}

\begin{proof}
The proof of (i) is similar to the proof of Lemma \ref{Kgen}. We have
(for $\mu,\lambda$ varying over $W\pi$)
\begin{equation}
\begin{aligned}
p_i(K_i) &= \prod_{\mu} (K_i - v_i^{\mu_i}) \\ 
&= \prod_{\mu} \left(\sum_\lambda v_i^{\lambda_i} \i_\lambda - v_i^{\mu_i}
\sum_\lambda \i_\lambda \right) \\
&= \prod_{\mu}\sum_\lambda (v_i^{\lambda_i} - v_i^{\mu_i})\i_\lambda \\
&= \sum_\lambda\prod_{\mu} (v_i^{\lambda_i} - v_i^{\mu_i})\i_\lambda \\
&= 0.
\end{aligned}
\end{equation}
This proves part (i).

Part (ii) follows immediately from the defining relations
\ref{S:relations}(c)--\ref{S:relations}(f), since if we choose $m$
sufficiently large we must have $E_i^m \i_\lambda = 0$ and $F_i^m
\i_\lambda = 0$ for all $\lambda \in W\pi$.  (Note that $\lambda \pm
m\alpha_i \notin W\pi$ for $m \gg 0$ since $W\pi$ is a finite set.)
\end{proof}

\subsection{}
We may regard every $\S$-module as an $\U$-module by composition with
the quotient map $\U \to \S$. A simple $\S$-module must be simple as a
$\U$-module. 

\begin{prop}\label{ss}
$\S=\S(\pi)$ is a finite-dimensional semisimple algebra.
\end{prop}

\begin{proof}
We know that $\U$ has a triangular decomposition $\U = \U^- \U^0 \U^+$
where $\U^0$ (resp., $\U^-$, $\U^+$) is the subalgebra of $\U$
generated by the $K_i, K_i^{-1}$ (resp., $F_i$, $E_i$).  It follows
that $\S=\S(\pi)$ has a similar decomposition $\S = \S^- \S^0 \S^+$
where $\S^0$, $\S^-$, $\S^+$ are defined as the homomorphic images of
$\U^0$, $\U^-$, $\U^+$, respectively, under the quotient map $\U \to
\S$ from Proposition \ref{image}. By the preceding lemma, we see that
the algebras $\S^0$, $\S^-$, $\S^+$ are finite-dimensional.  It
follows that $\S$ is finite-dimensional.

Moreover, $\S$ is a $\U$-module (via the quotient map $\U \to \S$).
Finite-dimensional $\U$-modules are semisimple. Thus $\S$ is
semisimple as a $\U$-module.  Hence $\S$ must be semisimple as an
$\S$-module, which proves that $\S$ is a semisimple algebra.
\end{proof}

\begin{lem}\label{wtspace}
If $M$ is any finite-dimensional $\S$-module then the decomposition
$M = \oplus_{\lambda\in W\pi} \i_\lambda M$ is a weight space
decomposition of $M$ as $\U$-module. Moreover, $M$ has type $\bf1$
when regarded as $\U$-module.
\end{lem}

\begin{proof}
Let $M$ be a finite-dimensional $\S$-module. As a $\U$-module, we have
the weight space decomposition $M = \oplus_{\mu,\sigma}
M^\mu_\sigma$ (see \ref{modules}).  On the other hand, in
$\S$ we have the equality $1=\sum_{\lambda\in W\pi} \i_\lambda$,
which implies that $M = \oplus_{\lambda\in W\pi} \i_\lambda M$. 
Clearly we have inclusions $\i_\lambda M \subset M^\lambda_{\bf1}$
for all $\lambda \in W\pi$, since 
\begin{equation}
 K_i \i_\lambda = v_i^{\lambda_i} \i_\lambda  \qquad (i = 1, \dots, n).
\end{equation}
Thus $M$ has type $\bf1$ (viewed as a $\U$-module). 

Moreover, if $m \in M^\lambda_{\bf1}$ then by definition we have $K_i
m = v_i^{\lambda_i} m$ for all $i$, so by multiplication by $\i_\mu$ we
obtain
\begin{equation}
v_i^{\mu_i} \i_\mu m = v_i^{\lambda_i} \i_\mu m \qquad(i =1, \dots, n).
\end{equation}
It follows that $\i_\mu m = 0$ for any $\mu \ne \lambda$. Thus $m =
1\cdot m = \sum_{\mu\in W\pi} \i_\mu m = \i_\lambda m$. This proves
that $m \in \i_\lambda M$, which establishes the inclusion
$M^\lambda_{\bf1} \subset \i_\lambda M$. Combining this with the
opposite inclusion from the previous paragraph, we conclude that
$M^\lambda_{\bf1} = \i_\lambda M$ for all $\lambda \in W\pi$.

Thus $M = \oplus_{\lambda\in W\pi} \i_\lambda M$ is the weight space
decomposition of $M$.
\end{proof}

\begin{prop}\label{iso}
The set $\{ \Lambda_\lambda \mid \lambda \in \pi \}$ is the set of
isomorphism classes of simple $\S$-modules., and $\dim \S(\pi) =
\sum_{\lambda \in \pi} (\dim \Lambda_\lambda)^2$.
\end{prop}

\begin{proof}
The simple $\S$-modules are simple $\U$-modules of type $\bf1$. Let
$\lambda \in X^+ - \pi$. If $\Lambda_\lambda$ was an $\S$-module, then
by Lemma \ref{wtspace} $\Lambda_\lambda$ would be a direct sum of the
weight spaces $\Lambda_\lambda^\mu$ as $\mu$ varies over $W\pi$. This
is a contradiction since $\Lambda_\lambda^\lambda \ne 0$ and $\lambda
\notin W\pi$.

On the other hand, for every $\lambda \in \pi$, $\Lambda_\lambda$
inherits a well-defined $\S$-module structure from its $\U$-module
structure, just by defining the action on the generators of $\S$ by
the obvious formulas. The first part of the proposition is proved.

The second part of the proposition follows immediately by standard
theory of finite-dimensional algebras.
\end{proof}

\subsection{} 
Let $\pi$ be an {\em arbitrary} subset of $X^+$. We could define an
algebra $\S(\pi)$ by the generators and relations given in
\sect\ref{Spi}. Then all of the results of this section, up to but not
including the preceding proposition, remain valid for that algebra,
since we do not use the assumption of saturation in any argument.
However, examining the proof of the preceding proposition, we see that
if $\pi$ is not saturated, then the only $\Lambda_\lambda$ which admit
an $\S(\pi)$-module structure are those of highest weight $\lambda$
belonging to $\pi'$, the largest saturated subset of $\pi$. In effect,
if $\pi$ is not saturated, the algebra $\S(\pi)$ collapses to
$\S(\pi')$. In particular, if $\pi$ has no saturated subset, then
$\S(\pi)$ is the zero algebra.  That is why we assumed, at the outset,
that $\pi$ is saturated.

\subsection{}
Following Donkin, we say that a $\U$-module {\em belongs to} $\pi$ if
all its composition factors are of the form $\Lambda_\lambda$ for
$\lambda \in \pi$.

\begin{cor}\label{GSA}
The algebra $\S(\pi)$ is isomorphic with the quotient algebra $\U/ I$
where $I$ is the ideal of $\U$ consisting of all elements of $\U$
which annihilate every $\U$-module belonging to $\pi$. Thus $\S(\pi)$
is a generalized $q$-Schur algebra in the sense of
\cite[\sect3.2]{Donkin:SA1}.
\end{cor}

\begin{proof}
Let $I$ be the kernel of the quotient mapping $\U \to \S$. From
Proposition \ref{iso} it is clear that $I$ annihilates every module
belonging to $\pi$.  On the other hand, if $x\in \U$ annihilates every
$\U$-module belonging to $\pi$ then the algebra $\U/(x)$ has at least
as many simple modules as does $\S = \U/I$, and thus there is a
quotient map $\U/(x) \to \U/I$, so $(x) \subset I$.
\end{proof}

\section{The algebra ${\bf\dot{\U}}/{\bf\dot{\U}}[P]$} 

\subsection{}
Set $P = X^+ - \pi$.  Then $P$ is a cofinite coideal in the poset
$X^+$ ($\lambda \le \mu$ for $\lambda \in P$, $\mu \in X^+$ implies
$\mu \in P$ and $X^+ - P$ is a finite set). Under precisely these
assumptions on $P$, Lusztig \cite[29.2]{Lusztig:book} showed that the
subspace ${\bf\dot{\U}}[P]$ of ${\bf\dot{\U}}$ spanned by $\sqcup_{\lambda \in
P} {\bf\dot{\B}}[\lambda]$ is a two-sided ideal of ${\bf\dot{\U}}$, and that the
corresponding quotient algebra ${\bf\dot{\U}}/ {\bf\dot{\U}}[P]$ is a
finite-dimensional semisimple algebra of dimension $\sum_{\lambda \in
\pi} (\dim \Lambda_\lambda)^2$.  The following result provides a
presentation by generators and relations for ${\bf\dot{\U}}/ {\bf\dot{\U}}[P]$.

\begin{thm} 
The algebra $\S(\pi)$ is isomorphic with Lusztig's algebra
${\bf\dot{\U}}/{\bf\dot{\U}}[P]$.
\end{thm}

\begin{proof}
First we show that ${\bf\dot{\U}}/{\bf\dot{\U}}[P]$ is a homomorphic image of
$\S(\pi)$. By Lusztig \cite[23.2.2(c)]{Lusztig:book} the algebra
${\bf\dot{\U}}$ is generated by elements of the form
$\divided{E_i}{m}1_\mu$, $\divided{F_i}{m}1_\mu$ for various $i = 1,
\dots, n$, $\mu \in X$, $m\ge 0$.  By \cite[29.2.1]{Lusztig:book} the
ideal ${\bf\dot{\U}}[P]$ contains the idempotents $1_\mu$ for $\mu \in WP =
X - W\pi$.  Write $\i_\mu$ for the image of $1_\mu$ under the quotient
map ${\bf\dot{\U}} \to {\bf\dot{\U}}/{\bf\dot{\U}}[P]$.  Thus $\i_\mu = 0$ for all
$\mu \in WP$, and thus in the quotient we have $\sum_{\lambda \in X}
\i_\lambda = \sum_{\lambda \in W\pi} \i_\lambda = 1$.

It follows that ${\bf\dot{\U}}/{\bf\dot{\U}}[P]$ is generated by all
$\i_\lambda$ ($\lambda \in W\pi$) and $E_i, F_i$ ($1\le i \le n$). The
elements $E_i, F_i$ are not elements of ${\bf\dot{\U}}$ but their images
make sense in the quotient since $E_i = \sum_{\lambda\in W\pi} E_i
\i_\lambda$, and similarly for $F_i$.

We need to verify that these generators satisfy the defining relations
of $\S(\pi)$. Relations \ref{S:relations}(a) follow from
\cite[23.1.1]{Lusztig:book}, \ref{S:relations}(b)--(f) follow from
\cite[23.1.3]{Lusztig:book}, and the quantized Serre relations
\ref{S:relations}(g), (h) hold on the $E_i$, $F_i$ since
they are, by Lusztig's definition of ${\bf\dot{\U}}$, images of the $E_i$,
$F_i$ in $\U$.  This proves that ${\bf\dot{\U}}/{\bf\dot{\U}}[P]$ is a
homomorphic image of $\S(\pi)$.  By dimension comparison, the two
algebras are isomorphic.  The proof is complete.
\end{proof}

\begin{cor}
The algebra $\S(\pi)$ has a canonical basis, formed by the non-zero
elements in the image of the canonical basis ${\bf\dot{\B}}$ of
${\bf\dot{\U}}$.
\end{cor}

\begin{proof}
This follows immediately from \cite[29.2.3]{Lusztig:book}.
\end{proof}

\subsection{} \label{Bdot}
Denote the canonical basis of $\S(\pi)$ by ${\bf\dot{\B}}(\pi)$.  It
may be regarded naturally as a subset of ${\bf\dot{\B}}$. It is the
disjoint union of the various two-sided cells ${\bf\dot{\B}}[\lambda]$ for
$\lambda \in \pi$.

\section{Change of base ring}\label{basechange}

\subsection{}
Define $\lsub{\S}{\A}$ to be the $\A$-subalgebra of $\S$ generated by
the elements $\i_\lambda$ ($\lambda\in W\pi$) and $\divided{E_i}{m} =
E_i^m/([m]_i^!)$, $\divided{F_i}{m} = F_i^m/([m]_i^!)$ ($1\le i \le
n$, $m\ge 0$).  Then $\lsub{\S}{\A}$ is the image of $\lsub{\U}{\A}$
under the quotient map (see Proposition \ref{image}) $\U \to \S$; it
is also the image of $\lsub{{\bf\dot{\U}}}{\A}$ under the map ${\bf\dot{\U}} \to
{\bf\dot{\U}}/{\bf\dot{\U}}[P] \simeq \S$.  Since ${\bf\dot{\B}}$ is an $\A$-basis
for $\lsub{{\bf\dot{\U}}}{\A}$, it follows that ${\bf\dot{\B}}(\pi)$ is a basis
for $\lsub{\S}{\A}$ as an $\A$-module.

\begin{prop}  
$\lsub{\S}{\A}$ is a cellular algebra, with cell datum
$(\pi,M,\b,\iota)$.  Its canonical basis ${\bf\dot{\B}}(\pi)$ is the
cellular basis for this cell datum.
\end{prop}

\begin{proof}
The cellularity follows immediately from the refined Peter-Weyl
theorem of Lusztig (see Theorem \ref{Peter-Weyl}).  It is immediate
that the anti-involution $\iota$ on ${\bf\dot{\U}}$ induces such a map on
the quotient $\S \simeq {\bf\dot{\U}}/{\bf\dot{\U}}[P]$.  The other properties
are easily checked.
\end{proof}

\subsection{}
Given any homomorphism $\A \to R$ to a commutative ring $R$ such that
$v$ maps to an invertible element of $R$, we may regard $R$ as an
$\A$-module.  We define $\lsub{\S}{R} = R \otimes_\A
\lsub{\S}{\A}$. Then $\lsub{\S}{R}$ is a cellular algebra with
essentially the same cell datum as $\lsub{\S}{\A}$; in particular, the
elements $1\otimes \b^\lambda_{S,T}$ for $\lambda\in \pi$, $S,T\in
M(\lambda)$ form a cellular basis over $R$.  Note that in case $R =
\Q(v)$ with $\A \to \Q(v)$ given by the obvious embedding, we have an
isomorphism $\lsub{\S}{\Q(v)} \simeq \S$.

\begin{thm}
If $R$ is a field then $\lsub{\S}{R}$ is quasihereditary.
\end{thm}

\begin{proof}
We apply the theory of cellular algebras from \cite{GL}.  According to
\cite[(3.10) Remark]{GL} it is enough to show that a certain bilinear
form $\phi_\lambda$ is nonzero for each $\lambda \in \pi$;
equivalently, it is enough to show that $\pi_0 = \pi$, where $\pi_0$
is the set of $\lambda \in \pi$ such that $\phi_\lambda \ne 0$.  By
\cite[(3.4) Theorem]{GL} the set $\pi_0$ is in bijective
correspondence with the set of isomorphism classes of simple
$\lsub{\S}{R}$-modules.

But one can easily construct a simple module for $\lsub{\S}{R}$ by
standard techniques, for each $\lambda \in \pi$. Set $\Delta(\lambda)
= R \otimes_\A \lsub{\Lambda_\lambda}{\A}$ where
$\lsub{\Lambda_\lambda}{\A} = \lsub{\U}{\A} v^+ = \lsub{{\bf\dot{\U}}}{\A}
v^+$ with $v^+ \ne 0$ a maximal vector ($\U^+v^+ = 0$) in
$\Lambda_\lambda$.  Then by standard arguments one sees that
$\Delta(\lambda)$ is a highest weight module; its unique simple
quotient is a highest weight module of highest weight $\lambda$.  One
obtains a simple module in this way for each $\lambda \in \pi$; these
simple modules are pairwise non-isomorphic.

This proves that $\pi_0 = \pi$; hence $\lsub{\S}{R}$ is quasihereditary. 
\end{proof}

\subsection{}
It is well-known that quasihereditary algebras have nice homological
properties; for instance, from the above it follows that (when $R$ is
a field) $\lsub{\S}{R}$ has finite global dimension and that its
matrix of Cartan invariants has determinant $1$ (see \cite[Theorem
1.1]{KX}).

\section{The classical case}\label{classical}

\subsection{} \label{u:relations}
Let $U$ be the associative algebra (with 1) on generators $e_i$,
$f_i$, $h_i$ ($1\le i \le n$) with relations
\begin{gather}
{}h_ih_j = h_jh_i, \quad e_if_j - f_je_i = \delta_{ij}h_i, \\
h_ie_j - e_jh_i = a_{ij}e_j, \quad h_if_j - f_jh_i = -a_{ij}f_j,\\ 
\sum_{s=0}^{1-a_{ij}} (-1)^s \binom{1-a_{ij}}{s}
e_i^{1-a_{ij}-s}e_je_i^s = 0 \quad (i \ne j) \\
\sum_{s=0}^{1-a_{ij}} (-1)^s \binom{1-a_{ij}}{s}
f_i^{1-a_{ij}-s}f_jf_i^s = 0 \quad (i \ne j). 
\end{gather}
Then $U$ is isomorphic with the universal enveloping algebra $U(\g)$
of the Lie algebra $\g$. The relations are the same as the relations
\ref{g:relations} defining $\g$, with the Lie bracket $[x,y]$ given by
$xy-yx$.

\subsection{} \label{s:relations}
We define an algebra $S(\pi)$ over $\Q$, depending only on the
saturated subset $\pi \subset X^+$ and the given Cartan matrix
$(a_{ij})$.  The algebra $S(\pi)$ is the associative algebra (with 1)
on generators $e_i$, $f_i$ ($1\le i \le n$), $\i_\lambda$ ($\lambda \in
W\pi$) with relations
\begin{gather}
\i_\lambda \i_\mu = \delta_{\lambda\mu} \i_\lambda, \quad
\sum_{\lambda\in W\pi} \i_\lambda = 1 \\
e_i f_j - f_j e_i = \delta_{ij} \sum_{\lambda\in W\pi}
\bil{\alpha_i^\vee}{\lambda} \i_\lambda \\ 
e_i \i_\lambda =
\begin{cases}
\i_{\lambda+\alpha_i} e_i &
   \text{if $\lambda+\alpha_i \in W\pi$}\\
0 & \text{otherwise}
\end{cases} \\
f_i \i_\lambda =
\begin{cases}
\i_{\lambda-\alpha_i} f_i &
   \text{if $\lambda-\alpha_i \in W\pi$}\\
0 & \text{otherwise}
\end{cases} \\
\i_\lambda e_i =
\begin{cases}
e_i \i_{\lambda-\alpha_i} &
   \text{if $\lambda-\alpha_i \in W\pi$}\\
0 & \text{otherwise}
\end{cases} \\
\i_\lambda f_i =
\begin{cases}
f_i \i_{\lambda+\alpha_i} &
   \text{if $\lambda+\alpha_i \in W\pi$}\\
0 & \text{otherwise}
\end{cases} \\
\sum_{s=0}^{1-a_{ij}} (-1)^s \binom{1-a_{ij}}{s}
e_i^{1-a_{ij}-s}e_je_i^s = 0 \quad (i \ne j) \\
\sum_{s=0}^{1-a_{ij}} (-1)^s \binom{1-a_{ij}}{s}
f_i^{1-a_{ij}-s}f_jf_i^s = 0 \quad (i \ne j).
\end{gather}
We obtained these relations by replacing $v$ by $1$ in
\ref{S:relations}.

\subsection{}
In the algebra $S = S(\pi)$ we define elements 
\begin{equation}
 h_i = \sum_{\lambda\in W\pi} \bil{\alpha_i^\vee}{\lambda}
\i_\lambda
\end{equation}
for $1\le i \le n$. From \ref{s:relations}(a) it follows that the
$e_i$, $f_i$, $h_i$ must satisfy relation \ref{u:relations}(a). In
particular, the $h_i$ commute with one another.

\begin{lem} Let $S^0$ be the subalgebra of $S$ generated by 
$h_1, \dots, h_n$. 

(i) The idempotents $\i_\lambda$ ($\lambda\in W\pi$) lie within $S^0$.

(ii) The $\i_\lambda$ ($\lambda\in W\pi$) form a basis for $S^0$.
\end{lem}

\begin{proof}\newcommand{\lp}{\lambda^\prime}
(Compare with the proof of Lemma \ref{Kgen}.)
Again we set $\lambda_i = \bil{\alpha_i^\vee}{\lambda}$ and write
$\Gamma(i,\lambda)=\{\mu\in W\pi \mid \mu_i = \lambda_i \}$. 
We set $J_i^\lambda = \prod_\mu (h_i - \mu_i)$ where the product is taken
over all $\mu \in W\pi - \Gamma(i,\lambda)$. 
We have equalities
\begin{equation}
\begin{aligned} 
J_i^\lambda &= \prod_\mu \left( \sum_{\lp \in W\pi} \lp_i \i_{\lp} -
\mu_i \sum_{\lp \in W\pi} \i_{\lp} \right) \\ &= \prod_\mu \left(
\sum_{\lp \in W\pi} (\lp_i - \mu_i) \i_{\lp} \right) \\ &= \sum_{\lp
\in W\pi} \prod_\mu (\lp_i - \mu_i) \i_{\lp} 
\end{aligned}
\end{equation}
where all
products are taken over $\mu \in W\pi - \Gamma(i,\lambda)$ and where
we have used the idempotent orthogonality relations
\ref{s:relations}(a) to interchange the sum and product.  Noting that
the product in the sum on the last line above vanishes for any $\lp
\in W\pi - \Gamma(i,\lambda)$, we obtain the expression
\begin{equation}
J_i^\lambda = \sum_{\lp
\in \Gamma(i,\lambda)} \prod_\mu (\lp_i - \mu_i) \i_{\lp}
\end{equation}
where the product in this sum is a nonzero {\em constant}, since
$\lp_i = \lambda_i$ for all $\lp \in \Gamma(i,\lambda)$.  This proves
that $J_i^\lambda$ is (up to a nonzero scalar) the sum of all
idempotents $\i_{\lp}$ for which $\lp_i = \lambda_i$. This property
holds for all $i$. Thus it follows that the product $J_1^\lambda
\cdots J_n^\lambda$ is, up to a nonzero scalar multiple, equal to
$\i_\lambda$, since $\i_\lambda$ is the unique idempotent appearing in
each of the sums in the product.  This proves part (i).

Part (ii) follows from part (i) by the same argument given in the
proof of \ref{Kgen}.
\end{proof}

\begin{prop}
The algebra $S=S(\pi)$ is a homomorphic image of $U$, via the map sending 
$e_i,f_i,h_i$ to the corresponding elements of $S$.
\end{prop}

\begin{proof}
Similar to the proof of Proposition \ref{image}.  From the lemma it
follows that $S(\pi)$ is generated by the elements $e_i,f_i,h_i$
($1\le i \le n$). These generators satisfy relations
\ref{u:relations}(a), as has been noted above. One easily verifies
that they satisfy relations \ref{u:relations}(b), by a calculation
similar to those given above.  They evidently satisfy
\ref{u:relations}(c), (d). The proof is complete.
\end{proof}

\subsection{}
Set $p_i(X) = \prod_{\mu\in W\pi} (X - \mu_i) \in \Q[X]$,
where $X$ is a formal indeterminate.

\begin{lem} \label{c:finiteness}
In the algebra $S(\pi)$ we have:
\par(i) The $h_i$ all satisfy the polynomial identity $p_i(h_i) = 0$
for $i = 1, \dots, n$. 
\par(ii) The $e_i$, $f_i$ are nilpotent.
\end{lem}

\begin{proof}
We have (for $\mu,\lambda$ varying over $W\pi$)
$$
\begin{aligned}
p_i(h_i) &= \prod_{\mu} (h_i - \mu_i) \\ 
&= \prod_{\mu} \left(\sum_\lambda \lambda_i \i_\lambda - \mu_i
\sum_\lambda \i_\lambda \right) \\
&= \prod_{\mu}\sum_\lambda (\lambda_i - \mu_i)\i_\lambda \\
&= \sum_\lambda\prod_{\mu} (\lambda_i - \mu_i)\i_\lambda \\
&= 0.
\end{aligned}
$$
This proves part (i).

Part (ii) follows immediately from the defining relations
\ref{s:relations}(c)--(f), since if we choose $m$ sufficiently large
we must have $e_i^m \i_\lambda = 0$ and $f_i^m \i_\lambda = 0$ for all
$\lambda \in W\pi$.
\end{proof}

\subsection{}
We denote by $L_\lambda$ the simple $U$-module of highest weight
$\lambda \in X^+$.  We may regard every $S$-module as a $U$-module by
composition with the quotient map $U \to S$. A simple $S$-module must
be simple as a $U$-module.

\begin{prop}\label{c:ss}
$S=S(\pi)$ is a finite-dimensional semisimple algebra.
\end{prop}

\begin{proof}
The proof is virtually identical to the proof of Proposition
\ref{ss}. 
\end{proof}

\begin{lem}\label{c:wtspace}
If $M$ is any finite-dimensional $S$-module then the decomposition
$M = \oplus_{\lambda\in W\pi} \i_\lambda M$ is a weight space
decomposition of $M$ as $U$-module. 
\end{lem}

\begin{proof}
The proof is similar to the proof of Lemma \ref{wtspace}. 
The details are left to the reader.
\end{proof}

\begin{prop}\label{c:iso}
The set $\{ L_\lambda \mid \lambda \in \pi \}$ is the set of
isomorphism classes of simple $S$-modules, and $\dim S(\pi) =
\sum_{\lambda \in \pi} (\dim L_\lambda)^2$.
\end{prop}

\begin{proof}
Again, the proof is similar to the proof given in the quantum case;
see \ref{iso}.
\end{proof}

\subsection{}
We say that a $U$-module belongs to $\pi$ if all its composition
factors are of the form $L_\lambda$ for $\lambda \in \pi$.

\begin{cor}\label{c:GSA}
The algebra $S(\pi)$ is isomorphic with the quotient algebra $U/ I$
where $I$ is the ideal of $U$ consisting of all elements of $U$ which
annihilate every $U$-module belonging to $\pi$. Thus $S(\pi)$ is a
generalized Schur algebra in the sense of \cite[\sect3.2]{Donkin:SA1}.
\end{cor}

\begin{proof}
The proof is the same as the proof of \ref{GSA}.
\end{proof}

\subsection{}
The cell datum for $\S(\pi)$ gives rise to a corresponding cell datum
for $S(\pi)$. Any element of $\lsub{\S}{\A}(\pi)$ can be written in
terms of an $\A$-linear combination of products of generators; such an
expression corresponds to a $\Z$-linear combination of products of
corresponding generators of $S(\pi)$, obtained from the $\A$-linear
combination by setting $v$ to $1$. In particular, the canonical basis
of $\S(\pi)$ determines a corresponding canonical basis of $S(\pi)$.

\subsection{}
One can define $\lsub{S}{\Z}(\pi)$ to be the subalgebra of $S(\pi)$
generated by the divided powers $e_i^m/(m!)$, $f_i^m/(m!)$ for various
$1 \le i \le n$, $m\ge 0$.  (It is the image of the Kostant $\Z$-form
$\lsub{U}{\Z}$ of $U$ under the quotient map $U \to S(\pi)$.)  The
cell datum for $S(\pi)$ is a cell datum for $\lsub{S}{\Z}(\pi)$.

For any ring $R$ we set $\lsub{S}{R}(\pi) = R\otimes_\Z
\lsub{S}{\Z}(\pi)$. Then it follows from Donkin's results that, when
$R$ is a field, $\lsub{S}{R}(\pi)$ is quasihereditary. Alternatively,
one can see this directly as a consequence of the cellular structure,
as we did above in the quantum case. (The argument is the same.)

\section{Examples} \label{examples}

\subsection{}\label{ex:tensor}
One natural and interesting class of examples can be constructed as
follows, for an arbitrary Cartan matrix. Fix a finite-dimensional
representation $V$ of $\U$ (resp., $U$) and for each $d\ge 0$ take
$\pi=\pi(d)$ to be the set of dominant weights occurring in
$V^{\otimes d}$. The resulting algebra $\lsub{\S}{R}(\pi)$ (resp.,
$\lsub{S}{R}(\pi)$) has connections with classical invariant theory.

A canonical special case of this construction is obtained by taking
$V$ to be the simple non-trivial module of {\em smallest} dimension.
In types $A$, $B$, $C$, $D$ this $V$ is the natural module (see
\cite[5A.1, 5A.2]{Jantzen:qbook} for an explicit contruction of this
module in the quantum case).  Denote the resulting algebra by
$\lsub{\S}{R}(d)$ (resp., $\lsub{S}{R}(d)$). It is defined for every
indecomposable Cartan matrix.

\subsection{}\label{ex:A}
In type $A_{n-1}$ the algebras $\lsub{\S}{R}(d)$, $\lsub{S}{R}(d)$
have been extensively studied; this is the motivating example.  The
algebra $\lsub{\S}{R}(d)$ is isomorphic with the $q$-Schur algebra
constructed by Dipper and James \cite{DJ1, DJ2}; the same algebra was
also considered by Jimbo \cite{Jimbo} in the generic case.  The
algebra $\lsub{S}{R}(d)$ is isomorphic with the classical Schur
algebra constructed in \cite{Green:book} in order to study the
polynomial representations of general linear groups (based on Schur's
dissertation).  The isomorphisms follows from \cite[1.4, 2.4]{DG:PSA}.
(Note that one needs to replace $q$ by $v=q^{1/2}$ in the Dipper-James
construction.)  The algebra $\lsub{\S}{\A}(d)$ was constructed and
studied from a geometric viewpoint in \cite{BLM}; it is the algebra
denoted ${\bf K}_d$ there.

There is a vast literature concerning these algebras. In particular,
$\lsub{\S}{\A}(d)$ can be constructed as the space of linear
endomorphisms of $\lsub{V}{\A}^{\otimes d}$ commuting with a natural
action of the Hecke algebra in type $A$; see \cite{Jimbo},
\cite{Du}. Similarly, $\lsub{S}{\Z}(d)$ is the space of linear
endomorphisms of $\lsub{V}{\Z}^{\otimes d}$ commuting with the action of the
symmetric group, acting by place permutation.  Here $\lsub{V}{\A}$
(resp., $\lsub{V}{\Z}$) are appropriate lattices in $V$.

\subsection{}\label{ex:C}
In type $C$ the algebra $\lsub{\S}{R}(d)$ is isomorphic with the
symplectic $q$-Schur algebra studied in \cite{Oehms}.  Indeed, the two
algebras in question have ``the same'' representation theory and the
same dimension.  Similarly, the algebra $\lsub{S}{R}(d)$ (in type $C$)
is isomorphic with the symplectic Schur algebra studied in
\cite{Donkin:symposium}, \cite{Doty:PRAMSACT}.

\subsection{} \label{ex:another}
There is another construction, equally natural, which leads to a
different generalization of Schur algebras. Again we fix a module $V$
and consider $V^{\otimes d}$. The action of $\U$ (resp., $U$)
determines a representation $\U \to \End(V^{\otimes d})$ (resp., $U
\to \End(V^{\otimes d})$) and one can study the finite-dimensional
algebra (and its various specializations) obtained as the image of the
representation. These types of algebras were called ``enveloping'' in
Weyl's book on the classical groups; in types $A$ and $C$ with $V$ the
natural module they coincide with the generalized Schur algebras
considered in \ref{ex:A}, \ref{ex:C}.

However, if one takes $V$ to be the adjoint module for $\U$ in type
$A_1$, and takes $d=1$, then the image of $\U \to \End(V)$ is a
$9$-dimensional algebra, while the generalized $q$-Schur algebra
$\lsub{\S}{R}(\pi)$ is a $10$-dimensional algebra. (It is the
$q$-Schur algebra in degree $2$.)  Thus the construction of this
subsection is, in general, different from that considered in
\ref{ex:tensor}.  Precisely the same difficulty occurs in the
classical setting.

Algebras of the type constructed in this subsection were considered in
\cite{Doty:PRAMSACT} for types $B$, $D$ (in the classical setting); I
do not know whether they are generalized Schur algebras. To prove that
they are, one needs to show that for every dominant weight of
$V^{\otimes d}$ there is a simple direct summand of $V^{\otimes d}$ of
that highest weight, where $V$ is the natural module.

\section{Integral forms} \label{integral}

\subsection{} \label{S_A:relations}
Consider the algebra (associative with $1$) over $\A$ generated by the
symbols $\divided{E_i}{a}$, $\divided{F_i}{b}$ ($1 \le i \le n, a,b \ge
0$), $\i_\lambda$ ($\lambda \in W\pi$) and satisfying the relations
\begin{gather}
\i_\lambda \i_\mu = \delta_{\lambda\mu} \i_\lambda, \quad
\sum_{\lambda\in W\pi} \i_\lambda = 1 \\
\divided{E_i}{a}\divided{E_i}{b} = \divided{E_i}{a+b}, \quad 
\divided{F_i}{a}\divided{F_i}{b} = \divided{F_i}{a+b}, \quad 
\divided{E_i}{0}=\divided{F_i}{0} = 1 \\
\divided{E_i}{a}\divided{F_j}{b} = \divided{F_j}{b}\divided{E_i}{a}
\quad (i\ne j) \\ 
\divided{E_i}{a} \i_{-\lambda} \divided{F_i}{b} = 
\sum_{t\ge 0} \sqbinom{a+b-\bil{\alpha_i^\vee}{\lambda}}{t}_i 
\divided{F_i}{b-t} \i_{-\lambda+(a+b-t)\alpha_i} \divided{E_i}{a-t}\\ 
\divided{F_i}{b} \i_\lambda \divided{E_i}{a} = 
\sum_{t\ge 0} \sqbinom{a+b-\bil{\alpha_i^\vee}{\lambda}}{t}_i 
\divided{E_i}{a-t} \i_{\lambda-(a+b-t)\alpha_i} \divided{F_i}{b-t}\\ 
\divided{E_i}{a} \i_\lambda =
\begin{cases}
\i_{\lambda+a\alpha_i}\divided{E_i}{a}  &
   \text{if $\lambda+a\alpha_i \in W\pi$}\\
0 & \text{otherwise}
\end{cases}\\
\divided{F_i}{b} \i_\lambda =
\begin{cases}
\i_{\lambda-b\alpha_i} \divided{F_i}{b} &
   \text{if $\lambda-b\alpha_i \in W\pi$}\\
0 & \text{otherwise}
\end{cases}\\
\i_\lambda \divided{E_i}{a} =
\begin{cases}
\divided{E_i}{a} \i_{\lambda-a\alpha_i} &
   \text{if $\lambda-a\alpha_i \in W\pi$}\\
0 & \text{otherwise}
\end{cases}\\
\i_\lambda \divided{F_i}{b} =
\begin{cases}
\divided{F_i}{b} \i_{\lambda+b\alpha_i} &
   \text{if $\lambda+b\alpha_i \in W\pi$}\\
0 & \text{otherwise}
\end{cases} \\
\sum_{s=0}^{1-a_{ij}} (-1)^s 
\divided{E_i}{1-a_{ij}-s}\divided{E_j}{1} \divided{E_i}{s} = 0 
\quad (i\ne j) \\
\sum_{s=0}^{1-a_{ij}} (-1)^s 
\divided{F_i}{1-a_{ij}-s}\divided{F_j}{1} \divided{F_i}{s} = 0 
\quad (i\ne j).
\end{gather}
Denote this algebra by $\lsub{\T}{\A}(\pi)$. It depends only on the
Cartan matrix $(a_{ij})$ and the saturated set $\pi$.  It follows from
Theorem \ref{iso:thm} ahead that $\divided{E_i}{a} = \divided{F_i}{b}
= 0$ for $a,b$ sufficiently large, so $\lsub{\T}{\A}(\pi)$ is finitely
generated.

\subsection{}
We define $\T(\pi) = \Q(v)\otimes_\A \lsub{\T}{\A}(\pi)$, where
$\Q(v)$ is regarded as an $\A$-algebra via the natural embedding $\A
\to \Q(v)$. We identify $\lsub{\T}{\A}(\pi)$ with the $\A$-submodule
of $\T(\pi)$ spanned by all elements of the form $1\otimes X$, for
$X\in \lsub{\T}{\A}(\pi)$; in this way $\lsub{\T}{\A}(\pi)$ may be
regarded as an $\A$-subalgebra of $\T(\pi)$.  In the algebra $\T(\pi)$
we write $E_i$ for $\divided{E_i}{1}$, $F_i$ for $\divided{F_i}{1}$.

\begin{thm}\label{iso:thm}
(i) There is an isomorphism $\S(\pi) \stackrel{\approx}{\longrightarrow}
\T(\pi)$ of $\Q(v)$-algebras sending $E_i$ to $E_i$, $F_i$ to $F_i$, and
$\i_\lambda$ to $\i_\lambda$.

(ii) The above isomorphism carries $\lsub{\S}{\A}(\pi)$ onto
$\lsub{\T}{\A}(\pi)$.
\end{thm}

\begin{proof}
As a $\Q(v)$-algebra, $\T(\pi)$ is the algebra generated by all
$\divided{E_i}{a}$, $\divided{F_i}{b}$ ($1 \le i \le n, a,b \ge 0$),
$\i_\lambda$ ($\lambda \in W\pi$) and satisfying the relations
\ref{S_A:relations}.  By \ref{S_A:relations}(b) we have the equalities
$\divided{E_i}{a}=E_i^a/([a]_i^!)$, $\divided{F_i}{b}=E_i^b/([b]_i^!)$
in $\T(\pi)$.  Thus the elements $E_i$, $F_i$, $\i_\lambda$ generate
$\T(\pi)$ as a $\Q(v)$-algebra.  These generators satisfy relations
\ref{S:relations}. Indeed, all of those relations are immediate,
except for \ref{S:relations}(b), which requires a small calculation
that we leave to the reader.  Therefore we know by Proposition 3.2 and
Theorem 4.1 that $\T(\pi)$ is a homomorphic image of $\U$ and also of
${\bf\dot{\U}}$, so that all relations in those algebras also hold in
$\T(\pi)$. Thus relations \ref{S_A:relations} are consequences of the
relations \ref{S:relations}. This follows from Lusztig
\cite[23.1.3]{Lusztig:book} and \cite[1.4.3]{Lusztig:book}.  Thus we
see that $\T(\pi)$ is determined by relations \ref{S:relations}, and
so there is an isomorphism $\S(\pi) \to \T(\pi)$ of $\Q(v)$-algebras,
taking $E_i \to E_i$, $F_i \to F_i$, $\i_\lambda \to \i_\lambda$.
This proves part (i).

By definition, $\lsub{\S}{\A}(\pi)$ is the subalgebra generated by
all $\divided{E_i}{a}$, $\divided{F_i}{b}$ ($1 \le i \le n, a,b \ge
0$), $\i_\lambda$ ($\lambda \in W\pi$). These map onto the
corresponding defining generators for $\lsub{\T}{\A}(\pi)$.  This
proves part (ii).
\end{proof}

\subsection{} \label{S_Z:relations}
Consider the $\Z$-algebra generated by symbols $\divided{e_i}{a}$,
$\divided{f_i}{b}$ ($1 \le i \le n, a,b \ge 0$), $\i_\lambda$
($\lambda \in W\pi$) and subject to the relations
\begin{gather}
\i_\lambda \i_\mu = \delta_{\lambda\mu} \i_\lambda, \quad
\sum_{\lambda\in W\pi} \i_\lambda = 1 \\
\divided{e_i}{a}\divided{e_i}{b} = \divided{e_i}{a+b}, \quad 
\divided{f_i}{a}\divided{f_i}{b} = \divided{f_i}{a+b}, \quad 
\divided{e_i}{0} = \divided{f_i}{0} = 1 \\
\divided{e_i}{a}\divided{f_j}{b} = \divided{f_j}{b}\divided{e_i}{a}
\quad (i\ne j) \\ 
\divided{e_i}{a} \i_{-\lambda} \divided{f_i}{b} = 
\sum_{t\ge 0} \binom{a+b-\bil{\alpha_i^\vee}{\lambda}}{t}
\divided{f_i}{b-t} \i_{-\lambda+(a+b-t)\alpha_i} \divided{e_i}{a-t}\\ 
\divided{f_i}{b} \i_\lambda \divided{e_i}{a} = 
\sum_{t\ge 0} \binom{a+b-\bil{\alpha_i^\vee}{\lambda}}{t}
\divided{e_i}{a-t} \i_{\lambda-(a+b-t)\alpha_i} \divided{f_i}{b-t}\\ 
\divided{e_i}{a} \i_\lambda =
\begin{cases}
\i_{\lambda+a\alpha_i}\divided{e_i}{a}  &
   \text{if $\lambda+a\alpha_i \in W\pi$}\\
0 & \text{otherwise}
\end{cases}\\
\divided{f_i}{b} \i_\lambda =
\begin{cases}
\i_{\lambda-b\alpha_i} \divided{f_i}{b} &
   \text{if $\lambda-b\alpha_i \in W\pi$}\\
0 & \text{otherwise}
\end{cases}\\
\i_\lambda \divided{e_i}{a} =
\begin{cases}
\divided{e_i}{a} \i_{\lambda-a\alpha_i} &
   \text{if $\lambda-a\alpha_i \in W\pi$}\\
0 & \text{otherwise}
\end{cases}\\
\i_\lambda \divided{f_i}{b} =
\begin{cases}
\divided{f_i}{b} \i_{\lambda+b\alpha_i} &
   \text{if $\lambda+b\alpha_i \in W\pi$}\\
0 & \text{otherwise}
\end{cases}\\
\sum_{s=0}^{1-a_{ij}} (-1)^s 
\divided{e_i}{1-a_{ij}-s}\divided{e_j}{1} \divided{e_i}{s} = 0 
\quad (i\ne j)\\
\sum_{s=0}^{1-a_{ij}} (-1)^s 
\divided{f_i}{1-a_{ij}-s}\divided{f_j}{1} \divided{f_i}{s} = 0 
\quad (i\ne j).
\end{gather}
Denote this algebra by $\lsub{T}{\Z}(\pi)$. It depends only
on the Cartan matrix $(a_{ij})$ and the saturated set $\pi$.

\subsection{}
We define $T(\pi) = \Q \otimes_\Q \lsub{T}{\Z}(\pi)$, where
$\Q$ is regarded as a $\Z$-algebra via the natural embedding $\Z
\to \Q$. We identify $\lsub{T}{\Z}(\pi)$ with the $\Z$-submodule
of $T(\pi)$ spanned by all elements of the form $1\otimes X$, for
$X\in \lsub{T}{\Z}(\pi)$; in this way $\lsub{T}{\Z}(\pi)$ may be
regarded as a $\Z$-subalgebra of $T(\pi)$.  In the algebra $T(\pi)$
we write $e_i$ for $\divided{e_i}{1}$, $f_i$ for $\divided{f_i}{1}$.

\begin{thm}
(i) There is an isomorphism $S(\pi) \stackrel{\approx}{\longrightarrow}
T(\pi)$ of $\Q$-algebras sending $e_i$ to $e_i$, $f_i$ to $f_i$, and
$\i_\lambda$ to $\i_\lambda$.

(ii) The above isomorphism carries $\lsub{S}{\Z}(\pi)$ onto
$\lsub{T}{\Z}(\pi)$.
\end{thm}

\begin{proof}
The proof is similar to the proof of the corresponding result in the
quantum case, given above. We leave the details to the reader.
\end{proof}

\end{document}